\documentclass[11pt]{article}
\usepackage{amsmath, array, amstext, amsthm, amssymb, amsfonts, latexsym, amscd, mathtools}

\usepackage{authblk}
\usepackage{tikz}
\usetikzlibrary{shadings,patterns,decorations.pathreplacing}
\newcolumntype{C}{>{$}c<{$}}
\newtheorem{theorem}{Theorem}[section]
\newtheorem{corollary}[theorem]{Corollary}
\newtheorem{lemma}[theorem]{Lemma}
\newtheorem{proposition}[theorem]{Proposition}

\newtheorem{remark}[theorem]{Remark}

\title{\bf Extremal graphs with minimum number of connected subgraphs in a given family}
\author[1]{Dinesh Pandey \thanks{Corresponding author: dpandey@wlu.ca} }
\author[2]{Peruvemba Sundaram Ravi \thanks{pravi@wlu.ca}}
\affil[1,2] {Lazaridis School of Business and Economics,
             Wilfrid Laurier University, 75 University Avenue West, Waterloo, Ontario, N2L3C5, Canada.}
\date{}

\begin{document}
\maketitle
\begin{abstract}
The subgraph number of a vertex in a graph is defined as the number of connected subgraphs containing that vertex. The graph and its vertex which correspond to the minimum subgraph number among all graphs on $n$ vertices and $k$ cut vertices have been characterised. Further, using this characterisation, the graphs with the minimum number of connected subgraphs among all graphs on $n$ vertices and $k$ cut vertices, with girth at least $k$, have been obtained. This turns out to characterise the graphs with the minimum number of connected subgraphs among all graphs on $n$ vertices and $k$ cut vertices for $0 \leq k \leq 4$.\\



\noindent {\bf Conflicts of Interest:} The authors declare no conflict of interest.
\end{abstract}

  \section{Introduction}
Counting connected subgraphs of a graph is a challenging problem in combinatorics. However, when the graph is a tree, this count simplifies to the number of subtrees, which is relatively easy to determine. Sz\'ekely and Wang initiated the study of extremal problems concerning the number of subtrees in various classes of trees in their papers \cite{Szekely} and \cite{Szekely 1}. Since then, researchers have characterised extremal trees with the maximum or minimum number of subtrees, in many tree families, for example \cite{Kirk, Li, Xiao} and  \cite{Zhang}. In general graphs, the analog of subtrees is connected subgraphs. When the graph is not a tree, there is no standard method for counting these subgraphs. For some specific graphs, such counts can be obtained. In \cite{Wilf}, Wilf derived a recurrence relation for the number of connected labeled graphs on $n$ vertices. Using this, an expression for the number of connected subgraphs of the complete graph $K_n$ is presented in \cite{Pandey}. However, even for a slight modification of $K_n$, such as $K_5-e$ (where $e$ is any edge of $K_5$), counting its connected subgraphs becomes more intricate. In general, there is no routine method for counting the connected subgraphs of a graph. Nevertheless, if a graph contains a cut vertex $w$, the number of its connected subgraphs can be expressed in terms of some specific connected subgraphs containing $w$. This insight enables the derivation of bounds for the number of connected subgraphs in certain families of graphs. In \cite{Pandey}, the authors identified graphs with the minimum and maximum number of connected subgraphs among unicyclic graphs on $n$ vertices, and among graphs on $n$ vertices with a fixed number of pendant vertices. This paper contributes further to the study of extremal problems related to the number of connected subgraphs. Specifically, we focus on the family of connected graphs on $n$ vertices with a fixed number of cut vertices. We address the following problem.

\begin{itemize}
\item Among all graphs on $n$ vertices and $k$ cut vertices, which graphs have the minimum number of connected subgraphs, and what is that minimum?
\end{itemize}

The problem is solved partially in this paper. We have completely characterised the graphs for $0\leq k \leq 4$ and determined the minimum number. In addition, when $k\geq 5$ we have characterised the graphs in the family of graphs on $n$ vertices and $k$ cut vertices with girth at least $k$, and have obtained the minimum number.\\  

Throughout the paper, the graphs are simple, finite, connected and undirected. Let $G$ be a graph with vertex set $V(G)$ and the edge set $E(G)$. The edge joining the vertices $u$ and $v$ of $G$ is denoted by $uv$. A {\it cut vertex} in $G$ is a vertex whose removal makes the graph disconnected. A vertex of degree one in $G$ is called a {\it pendant vertex}. An edge containing a pendant vertex of $G$ is called a {\it pendant edge}. For two isomorphic graphs $G_1$ and $G_2$, we use the notation $G_1\cong G_2$. The path and the cycle on $n$ vertices are denoted by $P_n$ and $C_n$, respectively. The complete bipartite graph $K_{1,n-1}$ is called a {\it star} on $n$ vertices. A {\it block} in $G$ is a maximal $2$-connected subgraph of $G$. A {\it pendant block} of $G$ is a block containing exactly one cut vertex of $G$. Two blocks of $G$ are {\it adjacent} if they share a cut vertex.  For $u,v\in V(G)$, the distance $d_G(u,v)$ or $d(u,v)$ is the number of edges in a shortest path joining $u$ and $v$. The number of connected subgraphs of $G$ is denoted by $F(G)$. For $v\in V(G)$, the number of connected subgraphs containing $v$ is called the {\it subgraph number} of $v$ and denoted by $f_G(v)$. In general, $F_G(v_1,v_2,\ldots v_{\ell})$ denotes the number of connected subgraphs of $G$ containing the vertices $v_1,v_2,\ldots, v_{\ell}$. $\mathfrak{C}_{n,k}$ denotes the family of  graphs on $n$ vertices with $k$ cut vertices. If $G$ is $2$-connected, then it has $0$ cut vertices. There is no graph on $n$ vertices with $n$ or $n-1$ cut vertices and the path $P_n$ is the only graph with $n-2$ cut vertices. So, if $G$ is a non path graph on $n$ vertices and $k$ cut vertices, then $0\leq k \leq n-3$. Therefore, we consider the family $\mathfrak{C}_{n,k}$ with $k\leq n-3$. To add further, if $G\in \mathfrak{C}_{n,k}$ and $G$ contains a block of size $r$, then $|V(G)|\geq k+r$, implying if a graph has $n-3$ cut vertices, then it can not contain a block of size more than $3$.\\
 
  
The paper is organised in the following way. In Section \ref{Prelim}, we define some terminology and results from the literature. The characterisation of the graph and its vertex which correspond to the minimum subgraph number among all graphs on $n$ vertices and $k$ cut vertices is discussed in Section \ref{vertex containing subgraphs}. In Section \ref{Subgraph index}, we explore the number of connected subgraphs of a graph in the family of graphs on $n$ vertices and $k$ cut vertices, and in the family of graphs on $n$ vertices and $k$ cut vertices with girth at least $k$. Finally, we briefly add concluding remarks, with possible directions for future research in Section \ref{Conclusion}. 
\section{Preliminaries}\label{Prelim}
An edge $e$ of $G$ is a connected subgraph of $G$, and for any $v\in V(G)$, there is a connected subgraph that contains both $v$ and $e$. Therefore, the following results hold.
\begin{lemma}\label{edge effect}
Let $G$ be a graph and $e\in E(G)$. Then
\begin{enumerate}
\item[$(i)$] for any $v\in V(G)$, $f_{G-e}(v)<f_G(v)$ and
\item[$(ii)$] $F(G-e)< F(G)$.
\end{enumerate}
\end{lemma}

Let $G$ be a graph and $w$ be a cut vertex of $G$. It is clear that there always exist two subgraphs $H_1$ and $H_2$, each with at least two vertices such that $G\cong H_1\cup H_2$ and $V(H_1)\cap V(H_2)=\{w\}$. Moreover, if $w$ is shared by more than two blocks, then there will be multiple such pairs of subgraphs. Such pairs of subgraphs are helpful in obtaining the number of connected subgraphs of $G$, and the subgraph number of a vertex in $G$. An appropriate pair can be chosen as required.

\begin{lemma}[\cite{Pandey}; \cite{Xiao}, Proposition 1]\label{Counting}
Let $G\cong G_1\cup G_2$ such that $V(G_1)\cap V(G_2)=\{w\}$. Then
\begin{enumerate}
\item[$(i)$] For any $v\in V(G_1)$, $f_G(v)=f_{G_1}(v)+f_{G_1}(v,w)(f_{G_2}(w)-1)$ and
\item[$(ii)$] $F(G)=F(G_1)+F(G_2)-1+(f_{G_1}(w)-1)(f_{G_2}(w)-1)$.
\end{enumerate}
\end{lemma}
\begin{corollary} \label{subgraphs containing w}
Let $G$ be a connected graph and $w$ be a cut vertex of $G$. If $G_1,G_2, \ldots, G_s$ are the subgraphs of $G$ such that $G\cong \cup_{i=1}^s{G_i}$ and $\cap_{i=1}^s {V(G_i)}=\{w\}$, then $f_G(w)=\prod_{i=1}^s f_{G_i}(w)$.
\end{corollary}
\begin{proof}
Proof follows by applying induction on $s$ and using Lemma \ref{Counting} $(i)$.
\end{proof}

If $G$ has a non pendant block $B$, then the number of connected subgraphs containing a vertex from $B$, and the number of connected subgraphs of $G$, can be expressed more precisely as in the following proposition. 

\begin{proposition}
Let $G\in \mathfrak{C}_{n,k}, k\geq 1$ and $B$ be a block in $G$ containing $s\leq k$ cut vertices. Let $w_1, w_2,\ldots, w_s$ be the cut vertices of $G$ in $B$. For $i=1,2,\ldots, s$, let $G_i$ be the maximal subgraph of $G$ such that $V(G_i)\cap V(B)=\{w_i\}$. Then 
\begin{enumerate}
\item[$(i)$] for any $v_0\in V(B)$,
\begin{align*}
f_G(v_0)&=f_B(v_0)+\sum_{i=1}^s{f_B(v_0,w_i)(f_{G_i}(w_i)-1)}+ \sum_{i<j}{f_B(v_0,w_i,w_j)(f_{G_i}(w_i)-1)(f_{G_j}(w_j)-1)}+\\
		&\:\:\:\;   \cdots\cdots\cdots + f_B(v_0,w_1,\ldots, w_s)(f_{G_1}(w_1)-1)(f_{G_2}(w_2)-1)\cdots (f_{G_s}(w_s)-1) \;\;and
\end{align*}
\item[$(ii)$] \begin{align*} F(G)&=F(B)+\sum_{i=1}^s{(F(G_i)-1)}+\sum_{i=1}^s{(f_B(w_i)-1)(f_{G_i}(w_i)-1)}+ \\
					      &\;\;\;\;\sum_{i<j}{f_B(w_i,w_j)(f_{G_i}(w_i)-1)(f_{G_j}(w_j)-1)}+\\
			           	      &\;\;\;\;\cdots\cdots\cdots + f_B(w_1,w_2\ldots, w_s)(f_{G_1}(w_1)-1)(f_{G_2}(w_2)-1)\cdots (f_{G_s}(w_s)-1).
\end{align*}

\end{enumerate}
\end{proposition}

\begin{proof}
\begin{enumerate}
\item [$(i)$]The connected subgraphs of $G$ containing $v_0$ are precisely
\begin{itemize}
\item the connected subgraphs containing $v_0$ lying entirely in $B$. The number of such subgraphs is $f_B(v_0)$.
\item the connected subgraphs containing $v_0$ and a vertex of $G_i$ different from $w_i$, for $1\leq i\leq s$. The number of such subgraphs is $\sum_{i=1}^s{f_B(v_0,w_i)(f_{G_i}(w_i)-1)}$.
\item the connected subgraphs containing $v_0$, a vertex from $G_i$ other than $w_i$, and a vertex from $G_j$ other than $w_j$, for $1\leq i< j\leq s$. The number of such subgraphs is $\sum_{i<j}{f_B(v_0,w_i,w_j)(f_{G_i}(w_i)-1)(f_{G_j}(w_j)-1)}$.\\
\vdots
\item the connected subgraphs containing $v_0$, and a vertex from each $G_i$ other than $w_i$, $1\leq i\leq s$. The number of such subgraphs is $f_B(v_0,w_1,\ldots, w_s)(f_{G_1}(w_1)-1)(f_{G_2}(w_2)-1)\cdots (f_{G_s}(w_s)-1)$.
\end{itemize}
$f_G(v_0)$ is the sum of the number of these connected subgraphs. 

\item[$(ii)$] The proof follows using a similar argument as in $(i)$.
\end{enumerate}
\end{proof}

We use many graphs in later sections which have paths, stars and cycles as connected subgraphs. Therefore, the number of connected subgraphs and the subgraph number of vertices of these graphs are important.  The following are known (see \cite{Pandey}).
\begin{enumerate}
\item[$(i)$] $F(P_n)= {n+1\choose 2}$ and if $P_n=v_1v_2\cdots v_n$, then $f_{P_n}(v_i)=i(n-i+1)$,
\item[$(ii)$] $F(K_{1,n-1})= 2^{n-1}+n-1$ and 
				\begin{equation*}
					f_{K_{i,n-1}}(v)= \begin{cases} 
									2^{n-1} &\mbox {if deg(v)=n-1},\\
									2^{n-2}+1 &\mbox{ if deg(v)=1}, 	 
				 				 \end{cases}
				\end{equation*}
\item[$(iii)$] $F(C_n)= n^2+1$ and for any $v\in V(C_n)$,  $f_{C_n}(v)=\frac{n^2+n+2}{2}$.
\end{enumerate}

\begin{figure}[h!]
\begin{center}
\begin{tikzpicture}[scale=.7]
\draw (0,0) circle [radius= 1cm];
\filldraw (1,0) circle [radius=.5mm];
\filldraw (2,0) circle [radius= .5mm];
\filldraw (3.3,0) circle [radius=.5mm];
\filldraw (4.3,0) circle [radius= .5mm];
\draw [dash pattern=on 1pt off 2pt] (2,0)--(3.3,0);
\draw (1,0)--(2,0);
\draw(3.3,0)--(4.3,0);
\draw(0,0) node {$C_g$};
\draw [decorate,decoration={brace,amplitude=5pt,mirror},xshift=2pt,yshift=0pt]
(0.85,-0.2) -- (4.3,-0.2) node [black,midway,yshift=-0.4cm]
{\footnotesize $P_{n-g+1}$};
\end{tikzpicture}
\caption{The graph $L_{n, g}$}\label{graph-Lnk}
\end{center}
\end{figure}
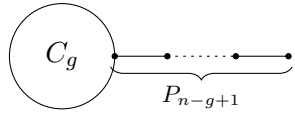

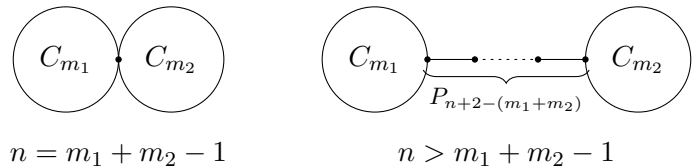
\begin{figure}[h!]
\begin{center}
\begin{tikzpicture}[scale =.7]
\draw (0,0)circle [radius= 1cm];
\draw (2,0) circle [radius= 1cm];
\draw  (1,-1.8) node {$n=m_1+m_2-1$};
\draw (0,0) node {$C_{m_1}$};
\draw (2,0) node {$C_{m_2}$};
\filldraw(1,0) circle [radius=.5mm];
\end{tikzpicture}
\hskip 1cm
\begin{tikzpicture}[scale=.7]
\draw (0,0)circle [radius= 1cm];
\draw (5,0) circle [radius= 1cm];
\filldraw (1,0)circle [radius= .5 mm];
\filldraw (1.9,0) circle [radius= .5 mm];
\filldraw (3.1,0)circle [radius= .5 mm];
\filldraw (4,0) circle [radius= .5 mm];
\draw [dash pattern= on 1pt off 2pt] (1.9,0)--(3.1,0);
\draw (1,0)--(1.9,0);
\draw(3.1,0)--(4,0);
\draw  (2.5,-1.8) node {$n>m_1+m_2-1$};
\draw (0,0) node {$C_{m_1}$};
\draw (5,0) node {$C_{m_2}$};
\draw [decorate,decoration={brace,amplitude=5pt,mirror},xshift=2pt,yshift=0pt]
(0.85,-0.2) -- (4,-0.2) node [black,midway,yshift=-0.4cm]
{\footnotesize $P_{n+2-(m_1+m_2)}$};
\end{tikzpicture}
\caption{The graphs $C_{m_1,m_2}^n$}\label{graph_Cmn}
\end{center}
\end{figure}

 By $L_{n,g}$, we denote the graph in $\mathfrak{C}_{n,n-g}$ obtained by identifying a pendant vertex of $P_{n-g+1}$ with a vertex of $C_g$ (see Figure \ref{graph-Lnk}). Note that $L_{n,n-g}$ has $g$ cut vertices. The graph $C_{m_1,m_2}^n$ is defined for $n\geq m_1+m_2-1$ and $m_1, m_2\geq 3$ in \cite{Pandey}. For $n\geq m_1+m_2$, $C_{m_1,m_2}^n$ is the graph obtained by identifying one pendant vertex of the path $P_{n+2-(m_1+m_2)}$ with a vertex of $C_{m_1}$ and the other pendant vertex of $P_{n+2-(m_1+m_2)}$ with a vertex of $C_{m_2}$. For $n=m_1+m_2-1$, $C_{m_1,m_2}^n$ is the graph obtained by identifying a vertex of $C_{m_1}$ with a vertex of $C_{m_2}$ (see Figure \ref{graph_Cmn}). Note that $C_{m_1, m_2}^n$ has $n+2-(m_1+m_2)$ cut vertices. The number of connected subgraphs of $L_{n,n-k}$ is

\begin{equation}\label {F(L(n,n-k))}
F(L_{n,n-k})=\frac{k}{2}(n^2+k^2- 2nk+n+3)+(n-k)^2+1. \;\;\;\; \mbox{ (see \cite{Pandey})}
\end{equation}

Let $v_0$ be the pendant vertex and $w$ be the vertex of degree $3$ in $L_{n,n-k}$. Then
\begin{equation}
f_{L_{n,n-k}}(v_0)=f_{P_k(v_0)}+f_{C_{n-k}}(w)=k+\frac{(n-k)^2+n-k+2}{2}=\frac{(n-k)^2+n+k+2}{2}.
\end{equation}

Consider the graph $C_{m_1,m_2}^n$, $n=m_1+m_2+k-2$ in which $w$ is the cut vertex lying on $C_{m_1}$. Then using Lemma \ref{Counting} ($ii$),
\begin{align}
F(C_{m_1,m_2}^n)&= F(C_{m_1})+F(L_{m_2+k-1,m_2})-1+(f_{C_{m_1}}(w)-1)(f_{L_{m_2+k-1,m_2}}(w)-1) \nonumber\\
			     &= m_1^2+m_2^2+1 + \frac{k-1}{2}(m_2^2+m_2+k+2)+(\frac{m_1^2+m_1+2}{2}-1)(\frac{m_2^2+m_2+2k}{2}-1)\nonumber\\
			     &=m_1^2+m_2^2+1 + \frac{k-1}{2}(m_2^2+m_2+k+2)+(\frac{m_1^2+m_1}{2})(\frac{m_2^2+m_2+2k-2}{2}) \nonumber\\
			     &=\frac{1}{4}(m_1^2m_2^2+m_1^2m_2+2m_1^2k+m_1m_2^2+m_1m_2+2m_1k+2m_1^2 \nonumber\\
			     &\;\;\;\;+2m_2^2+2m_2^2k +2m_2k+2k^2+2k-2m_1-2m_2).
			\label {Cm1m2n}
\end{align}

A graph is called minimally $2$-connected if it is $2$-connected and deleting any edge gives a graph which is not $2$-connected.
\begin{lemma}[\cite{Dirac}, Theorem 2] \label{min 2-connected}
A minimally $2$-connected graph with more than $3$ vertices is triangle free.
\end{lemma}

\section{Minimum subgraph number in $\mathfrak{C}_{n,k}$} \label{vertex containing subgraphs}
In this section, we obtain the graph and its vertex that correspond to the minimum subgraph number in $\mathfrak{C}_{n,k}$, which is used to prove our main results in Section \ref{Subgraph index}.
\begin{proposition}\label{no cut vertex}
Let $G\in \mathfrak{C}_{n,0}, n\geq 3$. Then for any $v\in V(G)$, $f_G(v)\geq \frac{n^2+n+2}{2}$ and equality holds if and only if $G\cong C_n$.
\end{proposition}
\begin{proof}
It is known that $f_{C_n}(u)=\frac{n^2+n+2}{2}$ for any $u\in V(C_n)$. Suppose $G\in \mathfrak{C}_{n,0}$, $G\ncong C_n$ and $v\in V(G)$. The vertex $v$, and $G$ itself are two connected subgraphs of $G$ containing $v$. Let $v'$ be a vertex of $G$ different from $v$. As $G$ is $2$-connected there are at least $2$ paths joining $v$ and $v'$. So there are at least $2(n-1)$ paths with $v$ as one pendant vertex. Now let $u,u'$ be two vertices of $G$ different from $v$. As $G$ is $2$-connected, $u,v,u'$ lies on a cycle and hence there is a path joining $u$ and $u'$ containing $v$. So there are at least ${n-1\choose 2}$ paths containing $v$ as an internal vertex. Now let $e$ be any edge of $G$. Then $G-e$ is a subgraph containing a cycle  and $v$. Thus $f_G(v)\geq 2+2(n-1)+{n-1\choose 2}+1>\frac{n^2+n+2}{2}$.
\end{proof}

Let $G$ be a graph and $v,v'\in V(G)$. The paths containing $v,v'$ in $G$ can be partitioned into the following four classes. $Q_G(vv')=:$ the set of paths in $G$ with pendant vertices $v$ and $v'$, $Q_G(vv'-)=:$ the set of paths containing $v, v'$ with one pendant vertex $v$ and the other pendant vertex different from $v'$, $Q_G(-vv')=:$ the set of paths containing $v,v'$ with one pendant vertex $v'$ and the other pendant vertex different from $v$ and $Q_G(-vv'-)=:$ the set of paths containing $v, v'$ with the pendant vertices different from $v$ and $v'$. 

\begin{lemma}\label{Cycle 2 vertices}
Let $u,v\in V(C_n)$. Then $\frac{n^2+2n+4}{4}\leq f_{C_n}(u,v)\leq \frac{n^2-n+4}{2}$. Moreover, the left equality holds iff $d(u,v)=\frac{n}{2}$ and the right equality holds iff $d(u,v)=1$. 
\end{lemma}

\begin{proof}
Let $u,v\in V(C_n)$. $C_n$ itself is a subgraph of $C_n$ containing $u,v$. Any other subgraph containing $u,v$ is a path and belongs to exactly one of the classes $Q_{C_n}(uv)$, $Q_{C_n}(uv-)$, $Q_{C_n}(-uv)$ or $Q_{C_n}(-uv-)$. Clearly, $|Q_{C_n}(uv)|=2$. If we fix $u$ as one pendant vertex and chose any vertex $x$ from $V(C_n)\setminus \{u,v\}$, we get a path containing $u,v$ with one pendant vertex as $u$ and the other as $x$. So $|Q_{C_n}(uv-)|=n-2$. Similarly $|Q_{C_n}(-uv)|=n-2$. Let $d(u,v)=d$. Then $C_n-\{u,v\}$ is a disjoint union of two paths $P_{d-1}$ and $P_{n-d-1}$. Let $V_1=V(P_{d-1})$ and $V_2=V(P_{n-d-1})$. Any path in $Q_{C_n}(-uv-)$ is a path with both pendant vertices either in $V_1$ or in $V_2$. So $|Q_{C_n}(-uv-)|={d-1\choose 2}+{n-d-1\choose 2}=\frac{1}{2}(n^2+2d^2-2nd-3n+4)$. This gives $f_{C_n}(u,v)=1+2+2(n-2)+\frac{1}{2}(n^2+2d^2-2nd-3n+4)=\frac{1}{2}(n^2+2d^2-2nd+n+2)$. As $1\leq d \leq \frac{n}{2}$ and $f_{C_n}(u,v)$ is strictly decreasing for $1\leq d\leq \frac{n}{2}$, it is maximum when $d=1$ and minimum when $d=\frac{n}{2}$. Hence the result follows.
\end{proof}
\begin{proposition}\label{u,v}
Let $G\in \mathfrak{C}_{n,0}$, $G\ncong C_n$ and $u, v\in V(G)$. Then $f_G(u,v)> f_{C_n}(u',v')$ for any $u',v'\in V(C_n)$.
\end{proposition}

\begin{proof}
If we show $f_G(u,v)>\frac{1}{2}(n^2-n+4)$, then the result follows from Lemma \ref{Cycle 2 vertices}. Let $u,v\in V(G)$. As $G$ is $2$-connected, $u,v$ lie on a cycle in $G$, implying $|Q_G(uv)|\geq 2$. Let $C_r$ be a smallest cycle in $G$ containing $u,v$. Then $C_r$ is a connected subgraph of $G$ containing $u,v$.  For any $x\in V(G)\setminus \{u,v\}$, there is at least one path containing $u,v$ with pendant vertices $u$ and $x$. If $x$ lies on $C_r$, then there is a path on $C_r$ containing $u,v$ with pendant vertices $u$ and $x$. If $x$ does not lie on $C_r$, take a vertex $x'$ on $C_r$ nearest to $x$. Then the path on $C_r$ containing $v,v$ with pendant vertices $u$ and $x'$ can be extended to get a path containing $u,v$ with pendant vertices $u$ and $x$. So $|Q_G(uv-)|\geq n-2$. Similarly $|Q_G(-uv)|\geq n-2$. Deleting any edge of $C_r$ we get a subgraph containing $u, v$ which is a graph containing a cycle. So, $f_G(u,v)\geq 2+1+2(n-2)+1=2n.$ \\

We show that there are at least ${n-2\choose 2}$ more connected subgraphs containing $u$ and $v$. Note that $n\geq 4$. Let $x,y\in V(G)\setminus \{u,v\}$, $x\neq y$ and let $z\in V(G)\setminus V(C_r)$ which is adjacent to  some vertex $z'\in V(C_r)$. Since $C_r$ is a smallest cycle containing $u,v$ and $G\ncong C_r$, such a $z$ exists. Further, at least one of $x$ or $y$ is different from $z$. Assume $z\neq x$. \\

\noindent{\bf Case I:} $x,y$ both lie on $C_r$\\
Take a path $P$ in $C_r$ containing $x,y,u,v,$ and $z'$ such that $x$ is a pendant vertex of the path. This path together with the edge $zz'$ is a connected subgraph of $G$ containing $u,v$ and not counted above.\\ 

\noindent{\bf Case II:} $x,y\in V(G)\setminus V(C_r)$\\
Let $P_{xx'}$ be a path from $x$ to $x'$ where $x'\in V(C_r)$ and $P_{y,y'}$ be a path from $y$ to $y'$ not containing $x$ (since $G$ is $2$-connected, such a path always exists), $y'\in V(C_r)$. Then $C_r\cup P_{x,x'}\cup P_{yy'}$ is a subgraph containing $u,v$ which is not counted above.\\

\noindent{\bf Case III:} $x\in V(C_r), y\in  V(G)\setminus V(C_r)$ or $y\in V(C_r), x\in  V(G)\setminus V(C_r)$\\
Without loss of generality assume that $x\in V(C_r)$ and $y\in V(G)\setminus V(C_r)$. $C_r\cup P_{yy'}$ is a subgraph containing $u,v$ not counted above.

Thus there are ${n-2\choose 2}$ more subgraphs containing $u,v$. So $f_G(u,v)\geq 2n +{n-2\choose 2}=\frac{1}{2}(n^2-n+6)>\frac{1}{2}(n^2-n+4)$. This completes the proof.
\end{proof}

\begin{lemma}\label{0,1 cut vertices}
Let $G\in \mathfrak{C}_{n,k}$, $k\in\{0,1\}$, $n\geq k+3$ and $G\ncong K_{1,3}$. Let $B$ be a block in $G$. Then for any $u,v\in V(B)$, $f_G(u,v)\geq 2(n-k)-1$. 
\end{lemma}
\begin{proof}
Let $G\in \mathfrak{C}_{n,0}$. Then by Proposition \ref{u,v}, $f_G(u,v)\geq f_{C_n}(u',v')$ for any $u',v'\in V(C_n)$ and by Lemma \ref{Cycle 2 vertices}, $f_{C_n}(u',v')\geq \frac{n^2+2n+4}{4}\geq 2(n-0)-1$, for $n\geq 4$. In $\mathfrak{C}_{3,0}$, $C_3$ is the only graph, and the result is trivially true.\\

Now suppose $G\in \mathfrak{C}_{n,1}$. If $G\cong K_{1,n-1}$, then $n\geq 5$ and $f_G(u,v)=2^{n-2}>2(n-1)-1$. Suppose $G\ncong K_{1,n-1}$ and it has $r$ blocks.  We use induction on $r$. Since $G$ has a cut vertex, $r\geq 2$. Let $w$ be the cut vertex of $G$.\\

\noindent Suppose $r=2$. Then $G\cong B\cup B'$ with $V(B)\cap V(B')=\{w\}$. Let $|V(B)|=n_1$ and $|V(B')|=n_2$. First suppose $n_2=2$. then $B\in \mathfrak{C}_{n_1,0}$ with $n_1\geq 3$, and $f_G(u,v)=f_B(u,v)+f_B(u,v,w)\geq 2n_1-1+2=2(n-1)+1>2(n-1)-1$. Now suppose $n_2\geq 3$. If $n_1=2$, then by Proposition \ref{no cut vertex}, $f_G(u,v)=f_{B'}(w)\geq \frac{(n-1)^2+n+1}{2}=\frac{n^2-n+2}{2}>2(n-1)-1$. If $n_1\geq 3$ then $f_G(u,v)=f_B(u,v)+f_B(u,v,w)(f_{B'}(w)-1)\geq 2(n_1-0)-1+\frac{n_2^2+n_2}{2}>2(n-1)-1$. So, the result is true for graphs with $2$ blocks in $\mathfrak{C}_{n,1}$.\\

\noindent Assume that the result is true for graphs having $r-1$ blocks.\\

\noindent Now, let $G$ have $r$ blocks. Choose a block $B'$ (different from $B$) of smallest size from the remaining $r-1$ blocks. Then there exists a subgraph $H$ of $G$ such that $G\cong H\cup B'$ with $V(H)\cap V(B')=\{w\}$. Let $|V(H)|=m_1$ and $|V(B')|=m_2$. 
\begin{align*}
f_G(u,v)&=f_H(u,v)+f_H(u,v,w)(f_{B'}(w)-1)\\
	    &> \begin{cases}
	                   2(m_1-1)-1+ 2   &\mbox {if $m_2=2$},\\
	                   2(m_1-1)-1+ 2m_2-1 &\mbox{if $m_2\geq 3$},
	                \end{cases}\\
	                &\geq 2(n-1)-1.
\end{align*}
This completes the induction step and hence the proof.
\end{proof}

special-pendant blocks or $s$-pendant blocks are introduced in \cite{Pandey-1}. An $s$-pendant block of $G$ is a pendant block which shares its cut vertex with exactly one non-pendant block. 

\begin{lemma}[\cite{Pandey-1}, Proposition 3.3]
Let $G\in \mathfrak{C}_{n,k}, k\geq 2$. Then $G$ has at least two vertex disjoint $s$-pendant blocks.
\end{lemma} 
 \noindent Although the statement of Proposition 3.3 in \cite{Pandey-1} does not explicitly mention vertex disjointness, it follows from the proof provided therein. \\

\begin{proposition}\label{min f(u,v)}
Let $G\in \mathfrak{C}_{n,k}$, $0\leq k\leq n-3$, $G\ncong K_{1,3}$. Let $B$ be a block in $G$. Then for any $u,v\in V(B)$, $f_G(u,v)\geq 2(n-k)-1$. 
\end{proposition}

\begin{proof}
We prove by induction on $k$. The base cases for $k=0,1$ hold by Lemma \ref {0,1 cut vertices}.\\

\noindent Assume that the result holds for graphs in $\mathfrak{C}_{n,s}$ for $s\leq k-1$.\\

\noindent Suppose $G\in \mathfrak{C}_{n,k}$, $k\geq 2$. Let $B'$ be an $s$-pendant block of $G$ different from $B$ ( in case $B$ is $s$-pendant choose $B'$ disjoint from $B$), containing the cut vertex $w'$. Further, take $H_2$ to be the subgraph of $G$ containing $B'$ and all the pendant blocks adjacent to $B'$, but no vertex other than $w'$ of a non-pendant block of $G$. Take $H_1\cong (G\setminus H_2)\cup \{w'\}$. Then $H_1$ contains $B$ and $G\cong H_1\cup H_2$ with $V(H_1)\cap V(H_2)=\{w'\}$.  Let $|V(H_1)|=n_1$ and $|V(H_2)|=n_2$.

If $H_2$ is a tree then $H_2\cong K_{1,n_2-1}$, and in that case, $f_{H_2}(w')=2^{n_2-1}>2(n_2-1)-1$. On the other hand, if $H_2$ contains a cycle, then by Lemma \ref{0,1 cut vertices}, $f_{H_2}(w')>f_{H_2}(w',w_0)\geq 2(n_2-1)-1$ (here $w_0$ can be any vertex of $H_2$ other than $w'$). Therefore in all cases, $f_{H_2}(w')> 2(n_2-1)-1$. \\

\noindent$H_1$ has $k-1$ cut vertices. First suppose $H_1\cong K_{1,3}$. In this case $k=2$ and $n=n_2+3$
\begin{align*}
f_G(u,v)&=f_{H_1}(u,v)+f_{H_1}(u,v,w')(f_{H_2}(w')-1)\\
	    &> 4+2(2(n_2-1)-2)\\
	    &=4n_2-4\\
	    &>2(n_2+3-2)-1 \;\;\;\;\mbox{for $n_2\geq 3.$}
\end{align*}
For $n_2=2$, $|V(G)|=5$ and it can be easily verified that $f_G(u,v)\geq 5$.\\

\noindent Now suppose $H_1\ncong K_{1,3}$. Clearly $n_1\geq k+1$. If $n_1 \geq k+2$, then by induction hypothesis, $f_{H_1}(u,v)\geq 2(n_1-k+1)-1$ and 
\begin{align*}
f_G(u,v)&=f_{H_1}(u,v)+f_{H_1}(u,v,w')(f_{H_2}(w')-1)\\
		&> 2(n_1-k+1)-1+2(n_2-1)-2\\
		&=2(n_1+n_2-1-k)-1\\
		&=2(n-k)-1.
\end{align*}  

If $n_1=k+1$, then $H_1\cong P_{n_1}$. If $n_1\geq 4$, then it is easy to check that  $f_{H_1}(u,v)=f_{P_{n_1}}(u,v)\geq n_1-1\geq 3$. and it follows that $f_G(u,v)\geq 2(n-k)-1$. If $n_1=3$, then $k=2$, and 
\begin{align*}
f_G(u,v)&=f_{P_3}(u,v)+f_{P_3}(u,v,w)(f_{H_2}(w)-1)\\
		&> 2+2(n_2-1)-2\\
		&= 2(n_2+2-2)-2\\
		&=2(n-2)-2,
\end{align*}  
which implies $f_G(u,v)\geq 2(n-2)-1$. This completes the proof.
\end{proof}

\begin{lemma}\label{Replace_cycle}
Let $G\in \mathfrak{C}_{n,1}$ and $B$ be a block of $G$, with $|V(B)|=n_1\geq 4$ and $B\ncong C_{n_1}$. Construct a new graph $G'$ from $G$ by replacing $B$ by $C_{n_1}$. Then for any $v\in V(G)$ there exists a $v'\in V(G')$ such that $f_{G'}(v')< f_G(v)$.
\end{lemma}
\begin{proof}
Let $w$ be the cut vertex of $G$ and $H$ be the subgraph induced by $V(G)\setminus V(B)\cup\{w\}$. Then $G\cong H\cup B$ with $V(H)\cap V(B)=\{w\}$ and $G'\cong H\cup C_{n_1}$ with $V(H)\cap V(C_{n_1})=\{w\}$. First suppose $v\in V(H)$. Then by Lemma \ref{Counting} $(i)$, 
\begin{align*}
f_G(v)&=f_H(v)+f_H(v,w)(f_{B}(w)-1)\\
	 &>f_H(v)+f_H(v,w)(f_{C_{n_1}}(w)-1) &\mbox{[using Proposition \ref{no cut vertex}]}\\
	 &=f_{G'}(v).
\end{align*}
Now suppose $v\in V(B)$. Chose any $v'\in V(C_{n_1})$ in $G'$. Then
\begin{align*}
f_G(v)&=f_{B}(v)+f_{B}(v,w)(f_H(w)-1)\\
	 &>f_{C_{n_1}}(v')+f_{C_{n_1}}(v',w)(f_H(w)-1) &\mbox{[using Proposition \ref{no cut vertex} and Proposition \ref{u,v}]}\\
	 &=f_{G'}(v').
\end{align*}
This completes the proof.
\end{proof}

\begin{lemma}\label{LC_replace}
Let $w_1$ be the cut vertex of $L_{n_1,n_1-1}$, $w_2$ be the cut vertex of $C_{m_1,m_2}^{n_1}$ where $n_1=m_1+m_2-1$, and $H$ be a connected graph containing a vertex $w$. Let $G_1$ be the graph obtained from $H$ and $L_{n_1,n_1-1}$ by identifying $w$ with $w_1$, $G_2$ be the graph obtained from $H$ and $C_{m_1,m_2}^{n_1}$ by identifying $w$ with $w_2$, and $G$ be the graph obtained from $H$ and $C_{n_1}$ by identifying $w$ with a vertex $w_3$ of $C_{n_1}$.  Then for any $v\in V(H)$, $f_G(v)< f_{G_i}(v)$,  $i=1,2$.  
\end{lemma}

\begin{proof}
Note that $G_1\cong H\cup L_{n_1,n_1-1}$ with $V(H)\cap V(L_{n_1,n_1-1})=\{w_1\}$, $G_2 \cong H\cup C_{m_1,m_2}^{n_1}$ with $V(H)\cap V(C_{m_1,m_2}^{n_1})=\{w_2\}$ and $G\cong H\cup C_{n_1}$ with $V(H)\cap V(C_{n_1})=\{w_3\}$. Using  Corollary \ref{subgraphs containing w}\\
$$f_{L_{n_1,n_1-1}}(w_1)=f_{C_{n_1-1}}(w_1)f_{P_2}(w_1)=n_1^2-n_1+2$$ and
\begin{align*}
f_{C_{m_1,m_2}^{n_1}}(w_2)&=f_{C_{m_1}}(w_2)f_{C_{m_2}}(w_2)\\
					    &=\frac{1}{4}(m_1^2m_2^2+m_1^2m_2+m_1m_2^2+m_1m_2+2m_1^2+2m_2^2+2m_1+2m_2+4).
\end{align*}
Further
\begin{align*}
f_{C_{n_1}}(w)&=\frac{n_1^2+n_1+2}{2}\\
		      &=\frac{1}{2}(m_1^2+m_2^2+2m_1m_2-m_1-m_2+2).
\end{align*}
So we get 
\begin{equation}\label{cycle-Ln1}
f_{C_{n_1}}(w)< f_{L_{n_1,n_1-1}}(w_1)
\end{equation}
and 
\begin{equation}\label{cycle-Cmn}
f_{C_{n_1}}(w)<f_{C_{m_1,m_2}^{n_1}}(w_2).
\end{equation}
 Now the result follows by counting $f_{G}(v), f_{G_1}(v)$ and $f_{G_2}(v)$ using Lemma \ref{Counting} $(i)$ and applying (\ref{cycle-Ln1}) and  (\ref{cycle-Cmn}).
\end{proof}

\begin{lemma}\label{min not cut}
Let $f_{G_0}(v_0)=\min\{f_G(v): G\in \mathfrak{C}_{n,k}, n\geq 4, v\in V(G)\}$. Then 
\begin{enumerate}
\item[$(i)$] $v_0$ is not a cut vertex of $G_0$ and
\item[$(ii)$] $G_0$ is triangle free.
\end{enumerate}
\end{lemma}

\begin{proof}\begin{enumerate}
\item[$(i)$] Let $v_0$ be a cut vertex of $G_0$. Then there exists subgraphs $G_1$ and $G_2$ of $G_0$ such that $v_0$ is a non cut vertex of $G_1$ and $G_0\cong G_1\cup G_2$ with $V(G_1)\cap V(G_2)=\{v_0\}$. Let $w$ be a non cut vertex of $G_1$ different from $v_0$. Construct a new graph $G_0'$ from $G_1$ and $G_2$ by identifying the vertex $w$ of $G_1$ with the vertex $v_0$ of $G_2$, i.e. $G_0'\cong G_1\cup G_2$ such that $V(G_1)\cap V(G_2)=\{w\}$. Then
\begin{align*}
f_{G_0}(v_0)&=f_{G_1}(v_0)f_{G_2}(v_0) & \mbox{[Corollary \ref{subgraphs containing w}]}\\
	     &=f_{G_1}(v_0)+f_{G_1}(v_0)(f_{G_2}(v_0)-1)\\
	     &>f_{G_1}(v_0)+f_{G_1}(v_0,w)(f_{G_2}(w)-1)\\
	     &=f_{G_0'}(v_0),
\end{align*}
which is a contradiction.
\item[$(ii)$] By Lemma  \ref{edge effect}$(i)$, blocks of $G_0$ are minimally $2$-connected and hence by Lemma \ref{min 2-connected} the blocks of size more than $3$ are triangle free. Let $B$ be a block of size $3$ in $G_0$. Then $B$ is isomorphic to a triangle, say $wxy$. As $n\geq 4$, $B$ must contain a cut vertex of $G_0$. Let $w$ be a cut vertex. If $B$ contains exactly one cut vertex $w$, then $G_0-xy \in \mathfrak{C}_{n,k}$ and by Lemma \ref{edge effect} $(i)$,  $f_{G_0-xy}(v_0)< f_{G_0}(v_0)$, a contradiction. If $B$ contains two cut vertices, say $w$ and $x$, then $G_0-wy \in \mathfrak{C}_{n,k}$ and $f_{G_0-wy}(v_0)< f_{G_0}(v_0)$, a contradiction. If all three vertices $w, x$ and $y$ are cut vertices in $G_0$, then for any $e\in \{wx, xy, yw \}$, $G_0-e \in \mathfrak{C}_{n,k}$ and $f_{G_0-e}(v_0)< f_{G_0}(v_0)$, a contradiction. 
\end{enumerate}
\end{proof}

Let $\mathfrak{T}_{n,k}$ be the set of all trees on $n$ vertices and $k$ cut vertices. Denote $\mathfrak{C}_{n,k}\setminus \mathfrak{T}_{n,k}$ by $\mathfrak{C'}_{n,k}$.

\begin{remark}\label{min not cut tree}
 Lemma \ref{min not cut} remains valid if the class $\mathfrak{C}_{n,k}$ is replaced by either $ \mathfrak{T}_{n,k}$, or by $\mathfrak{C'}_{n,k}$. 
\end{remark}

\begin{lemma}[\cite{Pandey},Lemma 4.6]\label{Lnk pendant_min}
Let $v_0$ be the pendant vertex and $v$ be a non pendant vertex of $L_{n,n-k}$. Then $f_{L_{n,n-k}}(v_0)<f_{L_{n,n-k}}(v)$.
\end{lemma}

\begin{theorem}\label{Ln,n-1}
Let $f_{G_0}(v_0)=\min\{f_G(v): G\in \mathfrak{C}_{n, 1}, n\geq 7, v\in V(G)\}$. Then $G_0\cong L_{n,n-1}$ and $v_0$ is the pendant vertex of $L_{n,n-1}$.
\end{theorem}
\begin{proof}

{\bf Claim I:} Every block of $G_0$ is either $K_2$ or a cycle on at least $4$ vertices. \\

Let $B$ be a block of $G_0$ on $n_1$ vertices which is neither $K_2$ nor $C_{n_1}$. By Lemma \ref{min not cut} $(ii)$, $n_1\geq 4$. Transform $G_0$ into a new graph $G_0'$ by replacing the block $B$ with the cycle $C_{n_1}$. Then by Lemma \ref {Replace_cycle}, there exists a $v_0'\in V(G_0')$ such that $f_{G_0'}(v_0')<f_{G_0}(v_0)$, a contradiction. So, $B\cong K_2$ or $B \cong C_{n_1}$.  This proves the claim I.\\

$G_0$ can not be isomorphic to $K_{1,n-1}$ because for $v\in V(K_{1,n-1})$ and $n\geq 7$, $f_{K_{1,n-1}}(v)\geq 2^{n-2}+1>\frac{n^2-n+4}{2}=f_{L_{n,n-1}}(v_0)$. So, $G_0$ has at least one cyclic block. \\

\noindent{\bf Claim II:} $G_0$ can not have more that two blocks.\\

Suppose $G_0$ has more than two blocks. Let $B_0$ be the block of $G_0$ containing $v_0$ and  $|V(B_0)|=n_0$. Let $w$ be the cut vertex of $G_0$ and $H$ be the subgraph induced by $(G_0\setminus B_0)\cup\{w\}$. Then $G_0\cong H\cup B_0$ with $V(H)\cap V(B_0)=\{w\}$. We now divide the possibilities in two cases.\\

\noindent{\bf Case I:} There is at least one cyclic block different from $B_0$\\
Let $B_1\cong C_{n_1}$ and $B_2$ with $|V(B_2)|=n_2$ be two blocks of $G_0$ different from $B_0$. Then $B_1\cup B_2$ is either isomorphic to $L_{n_1+n_2-1, n_1+n_2-2}$ or to $C_{n_1,n_2}^{n_1+n_2-1}$. Let $H'$ be the graph $(G_0\setminus B_1\cup B_2)\cup \{w\}$. Then $G_0$ is either isomorphic to $H'\cup L_{n_1+n_2-1, n_1+n_2-2}$ with $V(H')\cap V(L_{n_1+n_2-1, n_1+n_2-2})=\{w\}$ or to $H'\cup C_{n_1,n_2}^{n_1+n_2-1}$ with $V(H')\cap V(C_{n_1,n_2}^{n_1+n_2-1})=\{w\}$. Let $\bar{G_0}$ be the graph $H'\cup C_{n_1+n_2-1}$ such that $V(H')\cap V(C_{n_1+n_2-1})=\{w\}$. Then by Lemma \ref{LC_replace}, $f_{\bar{G_0}}(v_0)<f_{G_0}(v)$, a contradiction. \\

\noindent{\bf Case II} Except $B_0$ all other blocks are $K_2$\\
In this cases $G_0\cong C_{n_0}\cup K_{1,n-n_0}$ with $V(C_{n_0})\cap V(K_{1,n-n_0})=\{w\}$ and $n-n_0\geq 2$. Now
\begin{align*}
f_{G_0}(v_0)&=f_{C_{n_0}}(v_0)+f_{C_{n_0}}(v_0,w)(f_{K_{1,n-n_0}}(w)-1)\\
		   &\geq \frac{n_0^2+n_0+2}{2}+\left(\frac{n_0^2+2n_0+4}{4}\right)(2^{n-n_0}-1)\\
		   &\geq \frac{n_0^2+n_0+2}{2}+3\left(\frac{n_0^2+2n_0+4}{4}\right)\\
		   &= \frac{5n_0^2+8n_0+16}{4}.
\end{align*}
$n-n_0\leq 4$ otherwise by replacing $K_{1,n-n_0}$ by $C_{n-n_0+1}$, we get $C_{n_0,n-n_0+1}^n$ and $f_{C_{n_0,n-n_0+1}}(v_0)< f_{G_0}(v_0)$. So,
\begin{align*}
f_{L_{n,n-1}}(v_0)&=\frac{n^2-n+2}{2}\\
			   &\leq \frac{n_0^2+7n_0+14}{2}
\end{align*}
and 
\begin{equation*}
f_{G_0}(v_0)-f_{L_{n,n-1}}(v_0)\geq \frac{3n_0^2-6n_0-12}{4}>0,
\end{equation*}
Which is a contradiction. This proves claim II.\\

Therefore, $G_0$ has exactly two blocks with at least one a cycle. Therefore, either $G_0\cong L_{n,n-1}$ or $G_0\cong C_{m_1,m_2}^n$ for some $m_1,m_2\geq 4$ with $m_1+m_2-1=n$. Let $v$ be an arbitrary vertex in $C_{m_1,m_2}^n$. WLOG we may assume that $v$ lies on the cycle $C_{m_1}$. Then
\begin{align*}
f_{C_{m_1,m_2}^n}(v)&= f_{C_{m_1}}(v)+f_{C_{m_1}}(v,w)(f_{C_{m_2}}(w)-1)\\
				 &\geq \frac{m_1^2+m_1+2}{2}+(\frac{m_1^2+2m_1+4}{4})(\frac{m_2^2+m_2}{2})  &\mbox{\hskip -1.5cm[using Lemma \ref{Cycle 2 vertices}]}\\
				 &=\frac{1}{8}(4m_1^2+4m_2^2+2m_1m_2+m_1^2m_2^2+m_1^2m_2+2m_1m_2^2+4m_1+4m_2+8)\\
				 &>\frac{1}{8}(4m_1^2+4m_2^2+8m_1m_2-12m_1-12m_2+24)\\
				 &=f_{L_{n,n-1}}(v_0).
\end{align*}
This implies that $G_0\cong L_{n,n-1}$, and by Lemma \ref{Lnk pendant_min}, $v_0$ is the pendant vertex of $L_{n,n-1}$. 
\end{proof}
For $n\leq 6$, the number of graphs in $\mathfrak{C}_{n,1}$ is not very large. It can be easily verified that for $n\leq 6$, $f_{K_{1,n-1}}(v_0)=\min \{f_G(v): G\in \mathfrak{C}_{n,1}, v\in V(G)\}$, where $v_0$ is a pendant vertex of $K_{1,n-1}$. Furthermore, for $n\leq 5$, $f_{K_{1,n-1}}(v_0)$ uniquely attains the minimum, while for $n=6$, both $f_{K_{1,5}}(v_0)$ and $f_{L_{6,5}}(u_0)$, where $u_0$ is the pendant vertex of $L_{6,5}$, attain the minimum. \\  

Let $w$ be a cut vertex of $G$, and let $u$ and $v$ be two vertices that belong to different components of $G-w$. Then every connected subgraph containing both $u$ and $v$ must also contain $w$, whereas there are subgraphs that contain $u$ and $w$ but not $v$. Therefore, we have $f_G(u,w)>f_G(u,v)$ for such triplets of vertices $u,v,w$. Using this fact, the proof of the following lemma becomes straightforward. 

\begin{lemma}\label{Moving component}
Let $G\cong G_1\cup G_2 \cup G_3$ such that $V(G_1)\cap V(G_2)\cap V(G_3)=\{w\}$. Let $w'$ be a vertex of $G_2$ different from $w$. Construct a new graph $G^*$ from $G$ by moving $G_3$ to $w'$ from $w$ i.e. $G^*\cong G_1\cup G_2 \cup G_3$ such that $V(G_1\cup G_2)\cap V(G_3)=\{w'\}$ (see Figure \ref{fig}). Then for any $v\in V(G_1)$, $f_{G^*}(v)<f_G(v)$.   
\begin{figure}[h!]
\begin{center}

\begin{tikzpicture}[scale=0.7]
\filldraw (1,0) circle [radius=.3mm] (4.5,1)node[above]{$w'$} circle [radius=0.3mm];
\filldraw (3,0)node[left]{$w$} circle [radius=.3mm] (2,0)node[scale=.7] {$G_1$} (4,.7)node[scale=.7] {$G_2$} (4,-.7)node[scale=.7] {$G_3$} (2,-2.5)node{$G$} (2,-1)circle [radius=.3mm] node[below]{$v$};
\draw (1,0)..controls(2,2)..(3,0)..controls(2,-2)..(1,0);
\draw (3,0) to[out=-90,in=-155] (4.5,-1) to[out=-275,in=0] (3,0);
\draw (3,0) to[out=90,in=155] (4.5,1) to[out=275,in=0] (3,0);
\end{tikzpicture}
\hskip 1cm
\begin{tikzpicture}[scale=0.7]
\filldraw (1,0) circle [radius=.3mm] (4.5,1) circle [radius=0.3mm];
\filldraw (3,0)node[left]{$w$} circle [radius=.3mm] (2,0)node[scale=.7] {$G_1$}  (4,.7)node[scale=.7] {$G_2$} (5.25,1.5)node[scale=.7] {$G_3$} (2,-2.5)node{$G^*$} (4.35,1.4)node {$w'$} (2,-1)circle [radius=.3mm] node[below]{$v$};
\draw (1,0)..controls(2,2)..(3,0)..controls(2,-2)..(1,0);
\draw (4.5,1) to[out=90,in=155] (6,2) to[out=275,in=0] (4.5,1);
\draw (3,0) to[out=90,in=155] (4.5,1) to[out=275,in=0] (3,0);
\end{tikzpicture}
\end{center}
\caption{Moving a subgraph from one vertex to other}\label{fig}
\end{figure}
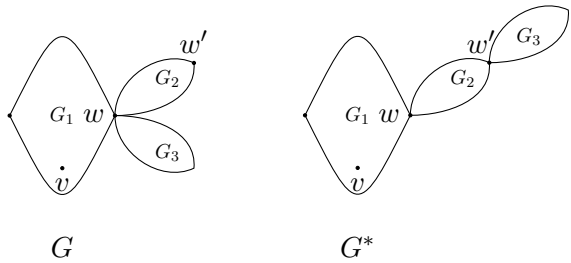

\end{lemma}

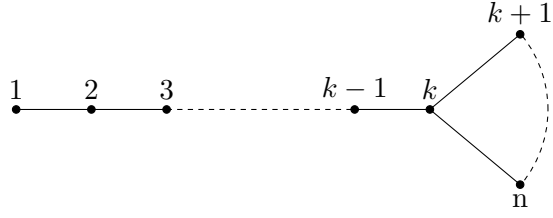
\begin{figure}[!h]
\begin{center}
\begin{tikzpicture}[scale =.5]
\filldraw (0,0)node[above]{1} circle [radius=1mm];
\filldraw (2,0)node[above]{2} circle [radius=1mm];
\filldraw (4,0)node[above]{3} circle [radius=1mm];
\filldraw (9,0)node[above]{$k-1$} circle [radius=1mm];
\filldraw (11,0)node[above]{$k$} circle [radius=1mm];
\filldraw (13.4,2)node[above]{$k+1$} circle [radius=1mm];
\filldraw (13.4,-2)node[below,outer sep=0pt]{n} circle [radius=1mm];
\draw(0,0)--(2,0)--(4,0);
\draw (9,0)--(11,0);
\draw (13.4,2)--(11,0)--(13.4,-2);
\draw [dash pattern= on 2pt off 2pt](4,0)--(9,0);
\clip (13.4,-2.6) rectangle (15,3);
\draw [dash pattern=on 2pt off 2pt] (11,0) circle [radius=3.14cm];
\end{tikzpicture}
\end{center}
\caption{The path-star tree $P_{k,n-k}$}\label{fig:2}
\end{figure}

The path-star tree (also known as broom or comet) on $n$ vertices with  $k-1$ cut vertices denoted by $P_{k,n-k}$ is the tree obtained by identifying one pendant vertex of the path $P_k$ with the center of the star $K_{1,n-k}$ (see Figure \ref{fig:2}). The following result is an immediate corollary of Lemma \ref{Moving component}.

\begin{corollary}\label{tree min vertex}
Let $v_0$ be the pendant vertex of the path-star tree $P_{k+1,n-k-1}$ lying on the path part. Then $f_{P_{k+1,n-k-1}}(v_0)=\min\{f_T(v): T\in \mathfrak{T}_{n,k}, v\in V(T)\}$. Furthermore, $f_{P_{k+1,n-k-1}}(v_0)=2^{n-k-1}+k$.
\end{corollary}
\begin{proof}
Let $f_{T_0}(v_0)=\min\{f_T(v): T\in \mathfrak{T}_{n,k}, v\in V(T)\}$. By Remark \ref{min not cut tree}, $v_0$ is a pendant vertex of $T_0$. By Lemma \ref{Moving component}, only one vertex of $T_0$ has degree more than $2$, and that vertex is a cut vertex farthest from $v_0$. This implies $T_0\cong P_{k+1,n-k-1}$ and $v_0$ is the pendant vertex lying on its path part.
\end{proof}

\begin{lemma}\label{sharing w}
Let $f_{G_0}(v_0)=\min\{f_G(v): G\in \mathfrak{C'}_{n,k}, n\geq 5, v\in V(G)\}$. Let $B$ be the block containing $v_0$ and $w$ be a cut vertex of $G_0$ in $B$. Then $B$ shares $w$ with at most four other blocks of $G_0$ and if it shares $w$ with two or more blocks, then all of them are pendant edges.
\end{lemma}

\begin{proof}
Suppose $B$ shares $w$ with two or more other blocks of $G_0$. Then there exist three subgraphs $G_1,G_2$ and $G_3$ of $G_0$ such that $G_0\cong G_1\cup G_2\cup G_3$ with $V(G_1)\cap V(G_2)\cap V(G_3)=\{w\}$ where $G_1$ contains $B$, and $w$ is a non-cut vertex in $G_2$ and $G_3$. Then by Lemma \ref{Moving component}, $G_2$ and $G_3$ are pendant blocks of $G_0$. Let $|V(G_2)|=n_2$ and $|V(G_3)|=n_3$. Assume that $n_2\geq n_3$. As $G_0$ is triangle free (by Lemma \ref{min not cut} $(ii)$), either $n_2=n_3=2$ or $n_2\geq 4$. \\

Suppose $n_2\geq 4$. Then by Proposition \ref{no cut vertex}, it follows that $G_2\cong C_{n_2}$ and $G_3$ is either isomorphic to $K_2$ or to $C_{n_3}$ for some $n_3\geq 4$. So, either  $G_0\cong G_1\cup L_{n_2+1,n_2}$ with $V(G_1)\cap V(L_{n_2+1,n_2})=\{w\}$ or $G_0\cong G_1\cup C_{n_2,n_3}^{n_2+n_3-1}$ with $V(G_1)\cap V(C_{n_2,n_3}^{n_2+n_3-1})=\{w\}$. Construct $G_0'$ from $G_0$ by replacing $G_2\cup G_3$ by $C_{n_2+n_3-1}$. Then $G_0'\in \mathfrak{C'}_{n,k}$. Since $f_{L_{n_2+1,n_2}}(w)>f_{C_{n_2+1}}(w)$ and $f_{C_{n_2,n_3}^{n_2+n_3-1}}(w)> f_{C_{n_2+n_3-1}}(w)$ (by (\ref{cycle-Ln1}) and (\ref{cycle-Cmn}), respectively), we get $f_{G_0}(v_0)>f_{G_0'}(v_0)$, which is a contradiction. So the only possibility is $G_2\cong K_2$ and $G_3\cong K_2$.

We proved that if $B$ shares $w$ with two or more other blocks, then all of them are pendant edges. Let $G_{11}$ be the maximal subgraph of $G_0$ containing $w$ as a non-cut vertex. Then $G_0\cong G_{11}\cup K_{1,\ell}$ with $V(G_{11})\cap V(K_{1,\ell})=\{w\}$. If $\ell\geq 5$, construct $G_0''$ from $G_0$ by replacing $K_{1,\ell}$ by $C_{\ell+1}$ such that $w$ remains a cut vertex shared by $G_{11}$ and $C_{\ell+1}$. Then it is easy to check that $f_{G_0''}(v_0)<f_{G_0}(v)$, a contradiction. This completes the proof.  
\end{proof}

\begin{theorem}\label{min subgraph number}
Let $G\in \mathfrak{C'}_{n,k}$, $1\leq k\leq n-3$. Then for any $v\in V(G)$, $f_G(v)\geq \frac{(n-k)^2+n+k+2}{2}$ and the equality holds if and only if $G\cong L_{n,n-k}$ and $v$ is the pendant vertex of $L_{n,n-k}$. 
\end{theorem}

\begin{proof}
Let $G\in \mathfrak{C'}_{n,k}$ and $v\in V(G)$. If $G\cong L_{n,n-k}$, then the result follows by Lemma \ref{Lnk pendant_min}. Now suppose $G\ncong L_{n,n-k}$. We use induction on $k$. Let $v_0$ be the pendant vertex of $L_{n,n-k}$ for $1\leq k\leq n-3$.\\

\noindent{\bf Base step:} Let $G \in \mathfrak{C'}_{n,1}$ and $v\in V(G)$. If $n\geq 7$, then by Theorem \ref {Ln,n-1}, $f_G(v)>  f_{L_{n,n-1}}(v_0)$. Also the result can be easily verified for $n=4,5,6$.\\

\noindent{\bf Induction hypothesis:} Let the result be true for $G\in \mathfrak{C'}_{n, \ell}$, $\ell\leq k- 1$.\\

\noindent{\bf Induction step:} Let $G\in \mathfrak{C'}_{n,k}$, $n\geq k+3$, and $v$ be an arbitrary vertex of $G$. Suppose $u_0$ is a vertex in $G$ such that $f_G(u_0)=\min\{f_G(x):x\in V(G)\}$. Then $f_G(v)\geq f_G(u_0)$. It is sufficient to show that $f_G(u_0)>f_{L_{n,n-k}}(v_0)$. Let $B$ be the block containing $u_0$ and $w$ be a cut vertex of $G$ in $B$. Let $G_1$ be the maximal subgraph of $G$ containing $B$ such that $w$ is a non cut vertex in $G_1$, and let $G_2$ be the maximal subgraph of $G$ containing $w$ and no other vertex of $G_1$. Thus $G\cong G_1\cup G_2$ with $V(G_1)\cap V(G_2)=\{w\}.$ Suppose $|V(G_1)|=n_1$ and $|V(G_2)|=n_2$. 

 
First suppose $|V(B)|\geq 3$. Then $G_1\in \mathfrak{C'}_{n_1,r}$ and $G_2\in \mathfrak{C}_{n_2,s}$ for some $r,s\geq 0$ such that either $k=r+s$ or $k=r+s+1$. As $G_2\in \mathfrak{C}_{n_2,s}$, $n_2\geq s+2$. Also since  $B$ contains a cycle, $n_1\geq r+3$. If $G_2\in \mathfrak{C'}_{n_2,s}$, then $n_2\geq s+3$ and $k=r+s+1$ (from Lemma \ref{sharing w}). So,
\begin{align*}
f_{G}(u_0)&= f_{G_1}(u_0)+f_{G_1}(u_0,w)(f_{G_2}(w)-1)\\
	  &\geq \frac{(n_1-r)^2+n_1+r+2}{2}+(2n_1-2r-1)( \frac{(n_2-s)^2+n_2+s}{2}) \\&\;\;\;\;\;\;\;\;\;\;\;\;\;\; \;\;\;\;\mbox{[using induction hypothesis and Proposition \ref{min f(u,v)} ]}\\
	  &=\frac{1}{2}((n_1-r)^2+(n_2-s)^2+2(n_1-r-1)(n_2-s)^2+n_1+r+2+(2n_1-2r-1)(n_2+s))\\	  
	  &>\frac{1}{2}((n_1-r-1+n_2-s)^2+ n_1+r+2+n_2+s)\\
	  &>\frac{1}{2}\{(n_1+n_2-1-(r+s+1))^2+ n_1+n_2+r+s+2\}\\
	  &=\frac{1}{2}((n-k)^2+n+k+2)\\
	  &=f_{L_{n,n-k}}(v_0).
\end{align*}
Now suppose $G_2\in \mathfrak{T}_{n_2,s}$. If $n_2\geq s+5$, then by Corollary \ref{tree min vertex}, $f_{G_2}(w)\geq 2^{n_2-s-1}+s \geq  \frac{(n_2-s)^2+n_2+s+2}{2}$. So as above we get, $f_{G}(u_0)> f_{L_{n,n-k}}(v_0)$. If $n_2\leq s+4$, then there are two possibilities for $G_2$. \\
{\bf Case I:} $w$ is a cut vertex of $G_2$\\
By Lemma \ref{sharing w}, $G_2\cong K_{1,n_2-1}$, $n_2\leq 5$. Also $r=k-1$ and 
\begin{align*}
f_{G}(u_0)&=f_{G_1}(u_0)+f_{G_1}(u_0,w)(f_{K_{1,n_2-1}}(w)-1)\\
	  &> \frac{(n_1-r)^2+n_1+r+2}{2}+(2(n_1-r)-1)( 2^{n_2-1}-1) \\&\;\;\;\;\;\;\;\;\;\;\;\;\;\; \;\;\;\;\mbox{[using induction hypothesis and Proposition \ref{min f(u,v)} ]}\\
	  &=\frac{(n_1+1-k)^2+n_1+k+1}{2}+(2(n_1-k)+1)(2^{n_2-1}-1)\\
	  &=\begin{cases}
	       \frac{(n-k)^2+n+k+2}{2}+27(n-k)-103 & \mbox{if $n_2=5$},\\
	       \frac{(n-k)^2+n+k+2}{2}+12(n-k)-35   & \mbox{if $n_2=4$},\\
	       \frac{(n-k)^2+n+k+2}{2}+5(n-k)-10     & \mbox{if $n_2=3$},\\
               \frac{(n-k)^2+n+k+2}{2}+2(n-k)-2       & \mbox{if $n_2=2$},
	  \end{cases}\\
	  &> f_{L_{n,n-k}}(v_0). 
\end{align*}
{\bf Case II:} $w$ is not a cut vertex of $G_2$ \\
In this case, by Corollary \ref{tree min vertex}, $G_2\cong P_{s+1,n_2-s-1}$ and $w$ is the pendant vertex lying on the path part of $P_{s+1,n_2-s-1}$. If $n_2-s-1\geq 4$, then by replacing $P_{s+1,n_2-s-1}$ by $L_{n_2,n_2-s}$ we get the graph $G_0'\cong G_1\cup L_{n_2,n_2-s}$ with $V(G_1)\cap V(L_{n_2,n_2-s})=\{w\}$, and it is easy to see that $f_{G_0}(u_0)>f_{G_0'}(u_0)$, which is a contradiction. Therefore $n_2-s-1\leq 3$. Now
\begin{align*}
f_G(u_0)&= f_{G_1}(u_0)+f_{G_1}(u_0,w)(f_{G_2}(w)-1)\\
	  &> \frac{(n_1-r)^2+n_1+r+2}{2}+(2n_1-2r-1)( 2^{n_2-s-1}+s-1) \\&\;\;\;\;\;\;\;\;\;\;\;\;\;\; \;\;\;\;\mbox{[using induction hypothesis and Proposition \ref{min f(u,v)} ]}\\
	    &=\begin{cases}
	              \frac{1}{2}((n_1-r)^2+4n_1s+5n_1-4rs-3r-2s)& \mbox{if $n_2-s-1=1$},\\
	              \frac{1}{2}((n_1-r)^2+4n_1s+13n_1-4rs-11r-2s-4)& \mbox{if $n_2-s-1=2$},\\
	              \frac{1}{2}((n_1-r)^2+4n_1s+29n_1-4rs-27r-2s-12)& \mbox{if $n_2-s-1=3$}.
	       \end{cases}\\
\end{align*}
In this case $k=r+s+1$, which gives
\begin{align*}
f_{L_{n,n-k}}(v_0)&=
	       \begin{cases}
	              \frac{1}{2}((n_1-r)^2+n_1+r+2s+4)& \mbox{if $n_2-s-1=1$},\\
	              \frac{1}{2}((n_1-r)^2+3n_1-r+2s+6)& \mbox{if $n_2-s-1=2$},\\
	              \frac{1}{2}((n_1-r)^2+5n_1-3r+2s+10)& \mbox{if $n_2-s-1=3$}.
	       \end{cases}\\
\end{align*}
Now 
\begin{align*}
f_{G}(u_0)-f_{L_{n,n-k}}(v_0)&=
	       \begin{cases}
	              \frac{1}{2}(4(n_1-r)s+4(n_1-r)-4s-4)& \mbox{if $n_2-s-1=1$},\\
	              \frac{1}{2}(4(n_1-r)s+10(n_1-r)-4s-10)& \mbox{if $n_2-s-1=2$},\\
	             \frac{1}{2}(4(n_1-r)s+24(n_1-r)-4s-22)& \mbox{if $n_2-s-1=3$},
	      \end{cases}\\ 
	      &>0. & \mbox{[as $n_1-r\geq 3$]}
\end{align*}

 Now suppose $V|(B)|=2$. By Remark \ref {min not cut tree}, $u_0$ is not a cut vertex. So, $u_0$ is a pendant vertex of $G_1$, implying $G_1\cong K_2$ and $G\cong K_2\cup G_2$ with $V(K_2)\cap V(G_2)=\{w\}$. By Lemma \ref{sharing w}, $G_2\in \mathfrak{C}_{n-1,k-1}$. Therefore,
 \begin{align*}
f_G(u_0)&=f_{K_2}(u_0)+f_{K_2}(u_0,w)(f_{G_2}(w)-1)\\
 	      &>2+\frac{(n-k)^2+n+k-2}{2} &\mbox{[using induction hypothesis]}\\
	      &=\frac{(n-k)^2+n+k+2}{2}\\
	      &=f_{L_{n,k}}(v_0).
 \end{align*}
 This completes the proof.
\end{proof}
 Note that
\begin{align*}
(2^{n-k-1}+k) \begin{cases}
			>\frac{(n-k)^2+n+k+2}{2}   &\mbox{if $k\leq n-6$},\\ 
			=\frac{(n-k)^2+n+k+2}{2}   &\mbox{if $k= n-5$},\\ 
	 		<\frac{(n-k)^2+n+k+2}{2}   &\mbox{if $k\geq n-4$}.
		       \end{cases}
\end{align*}
We conclude this section by the following result characterising the graphs and its vertex that correspond to the minimum subgraph number over $\mathfrak{C}_{n,k}$.
\begin{theorem}
Let $G\in \mathfrak{C}_{n,k}$, $1\leq k\leq n-3$ and $v\in V(G)$. 
\begin{itemize}
\item If $k\leq n-6$, then $f_G(v)\geq \frac{(n-k)^2+n+k+2}{2}$ and equality holds iff $G\cong L_{n,n-k}$ and $v$ is the pendant vertex of $L_{n,n-k}$.
\item If $k=n-5$, then $f_G(v)\geq 16+k$ and equality holds iff $G\cong P_{k+1,4}$ and $v$ is the pendant vertex of $P_{k+1,4}$ lying on its path part or $G\cong L_{n,5}$ and $v$ is the pendant vertex of $L_{n,5}$.
\item If $k\in\{n-3,n-4\}$, then $f_G(v)\geq 2^{n-k-1}+k$ and equality holds iff $G\cong P_{k+1,n-k-1}$ and $v$ is the pendant vertex of $P_{k+1,n-k-1}$ lying on its path part.
\end{itemize}
\end{theorem}

\begin{proof}
Let $u_0$ be the pendant vertex of $P_{k+1,n-k-1}$ lying on its path part and $v_0$ be the pendant vertex of $L_{n,n-k}$. Let $G\in \mathfrak{C}_{n,k}$ and $v\in V(G)$. Consider the following three cases.\\
\noindent {\bf Case I:} $k\leq n-6$. \\
If $G$ is a tree, then by Corollary \ref{tree min vertex}, $f_G(v)\geq 2^{n-k-1}+k> \frac{(n-k)^2+n+k+2}{2}$. If $G\in \mathfrak{C'}_{n,k}$, then by Theorem \ref{min subgraph number}, $f_G(v)\geq \frac{(n-k)^2+n+k+2}{2}$ and equality holds iff $G\cong L_{n,n-k}$ and $v=v_0$.\\

\noindent {\bf Case II:} $k= n-5$. \\
If $G$ is a tree, then by Corollary \ref{tree min vertex}, $f_G(v)\geq 16+k$ and equality holds iff $G\cong P_{k+1,n-k-1}$ and $v=u_0$. If $G\in \mathfrak{C'}_{n,k}$, then by Theorem \ref{min subgraph number}, $f_G(v)\geq 16+k$ and equality holds iff $G\cong L_{n,n-k}$ and $v=v_0$.\\

\noindent {\bf Case III:} $k\in \{n-4,n-3\}$.
If $G$ is a tree, then by Corollary \ref{tree min vertex}, $f_G(v) \geq 2^{n-k-1}+1$ and equality holds iff $G\cong P_{k+1,n-k-1}$ and $v=u_0$. If $G\in \mathfrak{C'}_{n,k}$, then by Theorem \ref{min subgraph number}, $f_G(v)\geq \frac{(n-k)^2+n+k+2}{2}>2^{n-k-1}+k$.
\end{proof}

\section{Graphs with the least number of connected subgraphs in $\mathfrak{C}_{n,k}(k)$}\label {Subgraph index}
In this section, we characterise the graphs with the minimum number of connected subgraphs among all graphs on $n$ vertices with $k$ cut vertices and girth at least $k$. Let $\mathfrak{C}_{n,k}(g)$ denote the subset of $\mathfrak{C}_{n,k}$ containing the graphs with girth at least $g$. The girth of a tree is assumed to be infinite and so $\mathfrak{T}_{n,k}\subsetneq \mathfrak{C}_{n,k}(g)$ for any finite $g$. We denote $\mathfrak{C}_{n,k}(g)\setminus \mathfrak{T}_{n,k}$ by $\mathfrak{C'}_{n,k}(g)$.  Observe that for $g=1,2,3$, we have $\mathfrak{C}_{n,k}(g)=\mathfrak{C}_{n,k}$ and $\mathfrak{C'}_{n,k}(g)=\mathfrak{C'}_{n,k}$.

\begin{lemma}[\cite{Pandey}, Proposition 4.11]\label{F0}
Let $G$ be a $2$-connected graph on $n\geq 3$ vertices. Then $F(G) \geq F(C_n)$ and the equality holds if and only if $G\cong C_n$.
\end{lemma}
Let $P_d$ ($d\geq 2$) be the path with pendant vertices $v_1$ and $v_d$. For $\ell, m\geq 1$, let $T(\ell,m,d)$ be the tree obtained by taking the path $P_d$ and adding $\ell$ vertices adjacent to $v_1$ and $m$ vertices adjacent to $v_d$. The tree $T(\ell,m,d)\in \mathfrak{C}_{l+m+d, d}$. Note that the connected subgraphs of a tree are precisely its subtrees.

\begin{lemma}[\cite{Li}, Theorem 1]\label{tree fixed pendant vertices}
Let $T$ be a tree on $n$ vertices and $k$ pendant vertices. Then $F(T)\geq F(T(\lfloor\frac{k}{2}\rfloor,\lceil\frac{k}{2}\rceil,n-k)$ and the equality holds iff $T\cong T(\lfloor\frac{k}{2}\rfloor,\lceil\frac{k}{2}\rceil,n-k)$. Furthermore,
    \begin{equation*}
F(T(\lfloor\frac{k}{2}\rfloor,\lceil\frac{k}{2}\rceil,n-k))=
\begin{cases}
(n-k-1)2^{\frac{k+2}{2}}+2^k+k+{n-k-1 \choose 2} & \textit{if k is even},\\ 
3(n-k-1)2^{\frac{k-1}{2}}+2^k+k+{n-k-1 \choose 2} & \textit{if k is odd}.
\end{cases}
\end{equation*}

\end{lemma}

\noindent As a tree on $n$ vertices with $k$ cut vertices has $n-k$ pendant vertices and vice versa, Lemma \ref{tree fixed pendant vertices} can be rephrased as the following.

\begin{proposition}\label{min on trees}
Let $T\in \mathfrak{T}_{n,k}$. Then $F(T)\geq F(T(\lfloor\frac{n-k}{2}\rfloor,\lceil\frac{n-k}{2}\rceil,k)$ and the equality holds iff $T\cong T(\lfloor\frac{n-k}{2}\rfloor,\lceil\frac{n-k}{2}\rceil,k)$. Furthermore,
\begin{equation*}
F(T(\lfloor\frac{n-k}{2}\rfloor,\lceil\frac{n-k}{2}\rceil,k)=\begin{cases}
										(k-1)2^\frac{n-k+2}{2}+2^{n-k}+n-k+{k-1 \choose 2} & \mbox{ if $n-k$ is even,}\\
										3(k-1)2^\frac{n-k-1}{2}+2^{n-k}+n-k+{k-1 \choose 2} & \mbox{ if $n-k$ is odd}.
									      \end{cases}
\end{equation*}
\end{proposition}


The proof of the following lemma proceeds by a similar argument to that used in the proof of Lemma \ref{min not cut}. 
\begin{lemma} \label{triangle free_F}
Let $F(G_0)=\min\{F(G): G\in \mathfrak{C}_{n,k}, k\geq 1, n\geq 4\}$. Then $G_0$ is triangle free. 
\end{lemma}

Next, we establish several results concerning the minimal graph with respect to the number of connected subgraphs in $\mathfrak{C'}_{n,k}(k)$ and $\mathfrak{C'}_{n,k}$, which will be instrumental in proving our main result. A key property of such minimal graphs is presented in the following lemma.

\begin{lemma}\label{pendant block is cycle}
Let $F(G_0)=\min\{F(G): G\in \mathfrak{C'}_{n,k}(k), k\geq 1, n\geq 4\}$. Then every pendant block of $G$ is either $K_2$ or a cycle. 
\end{lemma}
\begin{proof}
Let $B$ be a pendant block of $G_0$ which is not $K_2$. If $|V(B)|=3$, then $B$ is a triangle and the result follows. Suppose $|V(B)|=n_1\geq 4$ and $B\ncong C_{n_1}$. Let $w$ be the cut vertex of $G_0$ in $B$. Then there is a subgraph $H$ of $G_0$ such that $G_0\cong H\cup B$ with $V(H)\cap V(B)=\{w\}$. Construct a new graph $G_0'$ from $G_0$ by replacing $B$ with $C_{n_1}$, i.e.  $G_0'\cong H\cup C_{n_1}$ and $V(H)\cap V(C_{n_1})=\{w\}$. Then 
\begin{align*}
F(G_0)&=F(H)+F(B)-1+(f_H(w)-1)(f_B(w)-1)\\
	   &>F(H)+F(C_{n_1})-1+(f_H(w)-1)(f_{C_{n_1}}(w)-1) &\mbox{[Proposition \ref{no cut vertex} and Lemma \ref{F0}]}\\
	   &=F(G_0'),
\end{align*}
which is a contradiction. Hence $B\cong C_{n_1}$. This completes the proof.  
\end{proof}

\begin{lemma}\label{Comparision}
Let $C_{m_1,m_2}^n\in \mathfrak{C}_{n,k}(k)$. Then $F(L_{n,n-k})< F(C_{m_1,m_2}^n)$. 
\end{lemma}
\begin{proof}
Assume that $m_1\geq m_2$ in $C_{m_1,m_2}^n$. As $n=m_1+m_2+k-2$, (\ref{F(L(n,n-k))}) implies
\begin{align*}
F(L_{n,n-k})&=\frac{1}{4}(2m_1^2k+2m_2^2k+4m_1m_2k-6m_1k-6m_2k+2k^2+10k\\
	        &+ 4m_1^2+4m_2^2+8m_1m_2-16m_1-16m_2+20).
\end{align*}
Subtracting $F(L_{n,n-k})$ from $F(C_{m_1,m_2}^n)$ (computed in (\ref{Cm1m2n})), we get
\begin{align*}
F(C_{m_1,m_2}^n)- F(L_{n,n-k})&=\frac{1}{4} \left( m_1^2m_2^2+m_1^2m_2+m_1m_2^2-4m_1m_2k-7m_1m_2 +8m_1k+8m_2k \right.\\
						 &\;\;\;\; \left. -2m_1^2-2m_2^2-8k+14m_1+14m_2-20\right)\\
						 &=\frac{1}{4}\left( m_1m_2(m_1m_2-4k-7)+m_1^2(m_2-2)+m_2^2(m_1-2)\right.\\
						 &\;\;\;\; \left. +8(m_1+m_2-1)k+14(m_1+m_2)-20 \right)\\
						 &\geq \frac{1}{4}\left( m_1m_2(m_1m_2-4k-7)+40k+82 \right).
\end{align*}
First suppose $k\geq 6$. As $C_{m_1,m_2}^n\in \mathfrak{C}_{n,k}(k)$, $m_1m_2\geq 6k$, and so
\begin{align*}
F(C_{m_1,m_2}^n)- F(L_{n,n-k})&\geq m_1m_2(2k-7)\\
						 &>0.
\end{align*}
Now suppose $1\leq k\leq 5$. If $m_1m_2\geq 4k+7$, then clearly $F(C_{m_1,m_2}^n)- F(L_{n,n-k})>0$. Let $m_1m_2< 4k+7$. Then it can be easily verified that $m_1m_2-4k-7\geq -10$. Therefore
\begin{align*}
F(C_{m_1,m_2}^n)- F(L_{n,n-k}) &\geq \frac{1}{4}\left( m_1m_2(m_1m_2-4k-7)+40k+82 \right)\\
						  &> \frac{1}{4}(-10(4k+7)+40k+82)\\
						 & >0.
\end{align*}
\end{proof}

\begin{theorem}\label{F00}
Let $F(G_0)=\min \{F(G): G\in \mathfrak{C'}_{n,1}, n\geq 4\}$.  Then $G_0\cong L_{n,n-1}$. 
\end{theorem}
\begin{proof}
Let $B$ be a block of largest size in $G_0$. Then $B\cong C_{n_1}$ for some $3\leq n_1< n$. Suppose $G_0$ has more than two blocks. Let $B'$ be a block of size $n_2$ different from $B$ in $G_0$. By Lemma \ref{pendant block is cycle}, $B'$ is either $K_2$ or a cycle.  There exists a subgraph $H$ of $G_0$ such that $G_0\cong H\cup B\cup B'$ with $V(H)\cap V(B\cup B')=\{w\}$, where $w$ is the cut vertex of $G_0$. $B\cup B'$ is either isomorphic to $C_{n_1,n_2}^{n_1+n_2-1}$ or to $L_{n_1+1, n_1}$, and $w$ is the cut vertex of $B\cup B'$. Construct a new graph $G_0'\in \mathfrak{C'}_{n,1}$ from $G_0$ by replacing $B\cup B'$ with $C_{n_1+n_2-1}$, i.e. $G_0'\cong H\cup C_{n_1+n_2-1}$ with $V(H)\cap V(C_{n_1+n_2-1})=\{w\}$. 

It can be easily checked that $F(L_{n_1+1,n_1})>F(C_{n_1+1})$. Then using Lemma \ref{Comparision}, $F(C_{n_1,n_2}^{n_1+n_2-1})>F(C_{n_1+n_2-1})$. Also from (\ref{cycle-Ln1}) and \ref{cycle-Cmn}, we have  $f_{L_{n_1+1,n_1}}(w)>f_{C_{n_1+1}}(w)$  and $f_{C_{n_1,n_2}^{n_1+n_2-1}}(w)>f_{C_{n_1+n_2-1}}(w)$. Applying these inequalities to compare $F(G_0)$ and $F(G_0')$ using \ref{Counting} (ii), it  follows that $F(G_0)>F(G_0')$, which is a contradiction. Hence $G_0$ has exactly two blocks. Therefore either $G_0\cong L_{n,n-1}$ or $G_0\cong C_{n_1,n_2}^n$ for some $n_1,n_2\geq 3$ satisfying $n_1+n_2-1=n$. From Lemma \ref{Comparision}, $F(C_{n_1,n_2}^n)-F(L_{n,n-1})>0$. Hence $G_0\cong L_{n,n-1}$. This completes the proof.
\end{proof}

\begin{lemma}\label{cyclic pendant blocks}
Let $F(G_0)=\min\{F(G): G\in \mathfrak{C'}_{n,k}(k),2\leq k\leq n-3\}$. Then no cyclic pendant block of $G_0$ shares its cut vertex with another pendant block.
\end{lemma}
\begin{proof}
Let $B$ be a cyclic pendant block of $G_0$ and $w$ be the cut vertex of $G_0$ in $B$. Suppose $B'$ is another pendant block of $G_0$ such that $V(B)\cap V(B')=\{w\}$. So, $G_0\cong H\cup B\cup B'$ with $V(H)\cap V(B\cup B')=\{w\}$ for some subgraph $H$ of $G_0$. By Lemma \ref{pendant block is cycle}, either $B'\cong K_2$ or $B'\cong C_{n_2}$ for some $n_2\geq 3$. Let $B\cong C_{n_1}$, $n_1\geq 3$. Then  either $B\cup B'\cong L_{n_1+1,n_1}$ or $B\cup B'\cong C_{n_1,n_2}^{n_1+n_2-1}$, with $w$ as its cut vertex. Construct a new graph $G_0'$ from $G_0$ by replacing $B\cup B'$ by $C_{n_1+n_2-1}$ i.e. $G_0'\cong H\cup C_{n_1+n_2-1}$ with $V(H)\cap V(C_{n_1+n_2-1})=\{w\}$. Then $G_0'\in \mathfrak{C'}_{n,k}$ and as in the proof of Theorem \ref{F00}, it  can be showed that $F(G_0)> F(G_0')$, which is a contradiction.
\end{proof}

 Let $Q_{n,k}$ be the graph obtained by identifying a vertex of $C_{n-k-1}$ with the pendant vertex lying on the path part of $P_{k,2}$ i.e. $Q_{n,k}\cong C_{n-k-1}\cup P_{k,2}$ and $V(Q_{n,k})\cap V(P_{k,2})=\{u_0\}$, where $u_0$ is the pendant vertex lying on the path part of $P_{k,2}$. Note that if $G\in \mathfrak{C}_{n,k}(g)$, then $|V(G)|\geq k+g$.
 
\begin{theorem} \label{finite girth}
Let $G\in \mathfrak{C'}_{n,k}(k)$, $1\leq k\leq n-3$. Then $F(G)\geq F(L_{n,n-k})$. Moreover, If $n\neq 2k+1$, then the equality holds if and only if $G\cong L_{n,n-k}$, and for $n=2k+1$ the equality holds if and only if $G\in \{L_{2k+1,k+1}, Q_{2k+1,k}\}$.
\end{theorem}
\begin{proof}
We use induction on $k$. \\
\noindent{\bf Base step:} For $k=1$, the result follows from Theorem \ref{F00}.\\

\noindent{\bf Induction hypothesis:} Suppose the result is true for $ G\in \mathfrak{C'}_{n,g}(g), g\leq k-1$. \\

\noindent{\bf Induction step:} Let $G_0=\min\{F(G): G \in \mathfrak{C'}_{n,k}(k)\}$, $k\geq 2$. Let $B$ and $B'$ be two vertex disjoint $s$-pendant blocks of $G_0$ with $|V(B)|\leq |V(B')|$. Let $w$ be the cut vertex of $G_0$ in $B$. We first show that $B$ must be isomorphic to $K_2$.\\

Suppose $B\ncong K_2$. Then by Lemma \ref{pendant block is cycle}, $B\cong C_{n_1}$ for some $n_1\geq k$. By Lemma \ref{cyclic pendant blocks}, $B$ shares its cut vertex $w$ with exactly one other block, which is non-pendant. So, there exists a subgraph $H\in \mathfrak{C'}_{m,k-1}(k)$ such that $G_0\cong H\cup B$ with $V(H)\cap V(B)=\{w\}$, where $m=n-n_1+1$. By induction hypothesis, $F(H)\geq F(L_{m,m-k+1})$  and by Theorem \ref{min subgraph number}, $f_H(w)\geq F_{L_{m,m-k+1}}(w)$. So we get
\begin{align*}
F(G_0)&= F(C_{n_1})+F(H)-1+(f_{C_{n_1}}(w)-1)(f_H(w)-1)\\
	  &\geq F(C_{n_1})+F(L_{m,m-k+1})-1+(f_{C_{n_1}}(w)-1)(f_{L_{m,m-k+1}}(w)-1)\\
	  &= F(C_{n_1,m-k+1}^n)\\
	  &> F(L_{n,n-k}), & \mbox{[by Lemma \ref{Comparision}}] 
\end{align*}
which contradicts the minimality of $F(G_0)$. 

Therefore $B\cong K_2$. Let $H_1$ be the subgraph of $G_0$ containing $B$ and all the pendant blocks sharing the cut vertex $w$ but no vertex of a non-pendant block. As all the pendant blocks sharing $w$ are $s$-pendant, $H_1\cong K_{1,n_1-1}$ for some $n_1\geq 2$. Let $G_0\cong H_1\cup H_2$ with $V(H_1)\cap V(H_2)=\{w\}$, and let $|V(H_2)|=n_2$. As $G_0=\min\{F(G): G \in \mathfrak{C'}_{n,k}(k)\}$, $H_2$ must be $L_{n_2, n_2-k+1}$. If $n_1\geq 7$, then construct a new graph $G_0'$ from $G_0$ by replacing $H_1$ by $C_{n_1}$, then $G_0'\cong C_{n_1, n_2-k+1}^n$, and it can be checked that $F(G_0)>F(G_0')$, which contradicts the minimality of $F(G_0)$. So, $n_1\leq 6$. If $n_1=2$, then $G_0\cong L_{n,n-k}$ and the result follows. \\

Suppose $3\leq n_1\leq 6$. Let $w'$ be the cut vertex of $G_0$ lying on the cycle of $H_2 (\cong L_{n_2,n_2-k+1})$. Then $G_0\cong P_{k,n_1-1}\cup C_{n_2-k+1}$ with $V(P_{k,n_1-1})\cap V(C_{n_2-k+1})=\{w'\}$. Observing $G_0$ as this, we get

\begin{align*}
F(G_0)&=F(P_{k,n_1-1})+F(C_{n_2-k+1})-1+(f_{P_{k,n_1-1}}(w')-1)(f_{C_{n_2-k+1}}(w')-1)\\
           &= \frac{k(k-1)}{2}+2^{n_1-1}+n_1-1+(k-1)2^{n_1-1}+(n_2-k+1)^2\\
           &\;\;\;\;\;+(2^{n_1-1}+k-2)(\frac{(n_2-k+1)^2+n_2-k+1}{2})\\
           &=\frac{1}{2}((2^{n_1-1}+k)(n-n_1-k+2)^2+(2^{n_1-1}+k-2)(n-n_1-k+2)+k2^{n_1}+k(k-1)+2n_1-2)\\
           &=\begin{cases}
           	\frac{1}{2}((4+k)(n-k-1)^2+(k+2)(n-k-1)+k^2+7k+4) & \mbox{if $n_1=3$},\\
	        \frac{1}{2}((8+k)(n-k-2)^2+(k+6)(n-k-2)+k^2+15k+6) & \mbox{if $n_1=4$},\\
	        \frac{1}{2}((16+k)(n-k-3)^2+(k+14)(n-k-3)+k^2+31k+8) & \mbox{if $n_1=5$},\\
	        \frac{1}{2}((32+k)(n-k-4)^2+(k+30)(n-k-4)+k^2+63k+10) & \mbox{if $n_1=6$},
               \end{cases}
\end{align*}
and 
\begin{align*}
F(L_{n,n-k})&=  \frac{1}{2}((2+k)(n-k)^2+k(n-k)+k^2+3k+2).
\end{align*}
As $G_0\in \mathfrak{C'}_{n,k}(k)$ and  $n_1\geq 3$, $n\geq 2k+1$. Therefore we get,
\begin{align*}
F(G_0)-F(L_{n,n-k})&= \begin{cases}
           				 \frac{1}{2}((n-k)(2n-4k-6)+4k+4) & \mbox{if $n_1=3$},\\
	       				 \frac{1}{2}((n-k)(6n-10k-26)+14k+24) & \mbox{if $n_1=4$},\\
				         \frac{1}{2}((n-k)(14n-20k-82)+34k+108) & \mbox{if $n_1=5$},\\
	   			         \frac{1}{2}((n-k)(30n-38k-226)+72k+400) & \mbox{if $n_1=6$},
                                     \end{cases}\\
                               &\geq 0,
\end{align*}
and equality holds iff $n_1=3$ and $n=2k+1$. Therefore by minimality of $G_0$, $n_1=3$ and $n=2k+1$. In this case $G_0\cong Q_{2k+1,k}$ and $F(Q_{2k+1,k})=F(L_{2k+1,k+1})$.\\

Thus we conclude that If $n\neq 2k+1$, then $L_{n,n-k}$ uniquely minimizes the number of connected subgraphs, and if $n=2k+1$, then both $L_{2k+1,k+1}$ and $Q_{2k+1,k}$ attain the minimum. This completes the proof.
\end{proof}

\begin{theorem}\label{girth_k}
Let $G\in \mathfrak{C}_{n,k}(k)$, $n\geq 13$, $k\geq 1$. Then $F(G)\geq L_{n,n-k}$. Furthermore, If $n\neq  2k+1$, then the equality holds if and only if $G\cong L_{n,n-k}$, and for $n=2k+1$, the equality holds iff $G\in \{L_{2k+1, k+1}, Q_{2k+1,k}\}$.
\end{theorem}
\begin{proof}
Let $G\in \mathfrak{C}_{n,k}(k)$, $n\geq 13$.  Clearly $n\geq 2k$. So, $n\geq 13$ implies  $n-k\geq 7$. If $G\in \mathfrak{C'}_{n,k}(k)$, then the result follows from Theorem \ref{finite girth}. Suppose $G\in \mathfrak{C}_{n,k}(k)\setminus \mathfrak{C'}_{n,k}(k)$ i.e. $G$ is a tree. Then by Proposition \ref{min on trees}, $F(G)\geq F(T(\lfloor\frac{n-k}{2}\rfloor,\lceil\frac{n-k}{2}\rceil, k))$. We have,
\begin{align*}
&F(T(\lfloor\frac{n-k}{2}\rfloor,\lceil\frac{n-k}{2}\rceil,k)-F(L_{n,n-k})\\
&=\begin{cases}
	k(2^{\frac{n-k+2}{2}}-\frac{(n-k)^2+n-k}{2}-3)+2^{\frac{n-k+2}{2}}(2^{\frac{n-k-2}{2}}-1)-(n-k)(n-k-1), &\mbox {$n-k$ even},\\
	k(3^{\frac{n-k-1}{2}}-\frac{(n-k)^2+n-k}{2}-3)+2^{\frac{n-k-1}{2}}(2^{\frac{n-k+1}{2}}-3)-(n-k)(n-k-1), &\mbox {$n-k$ odd}.
     \end{cases}\\
     &>0\;\;\; \mbox{for $n-k\geq 7$}.
\end{align*}
Hence $F(G)>F(L_{n,n-k})$. This completes the proof.
\end{proof}

 \begin{table}[h!]
\begin{center}
\begin{tabular}{|c|c|c|c|c|c|c|}
\hline
$n$ & $\mathfrak{C}_{n,1}$ & $\mathfrak{C}_{n,2}$ & $\mathfrak{C}_{n,3}$ & $\mathfrak{C}_{n,4}(4)$ & $\mathfrak{C}_{n,5}(5)$ & $\mathfrak{C}_{n,6}(6)$\\
\hline
$12$ & $L_{12,11}; 190$ & $L_{12,10}; 216$ & $L_{12,9}; 226$ & $L_{12,8};223$ & $L_{12,7};210$ & $T(3,3,6);160$\\
\hline
$11$ & $L_{11,10};163$ & $L_{11,9};177$ &$L_{11,8};179$ & $L_{11,7}; 176$ & $T(3,3,5);140$ & $T(2,3,6);107$ \\
\hline 
$10$ & $L_{10,9};129$ & $L_{10,8};142$ & $L_{10,7};143$ &$T(3,3,4);121$ & $T(2,3,5);91$ & $T(2,2,6);70$\\
\hline
$9$ & $L_{9,8};103$ & $L_{9,7};111$ & $T(3,3,3);103 $ &$T(2,3,4);76$ & $T(2,2,5);58$ & $T(1,2,6);51$\\ 
\hline 
$8$ & $L_{8,7}; 80$ & $L_{8,6}; 84$ & $T(2,3,3); 62$ & $T(2,2,4);47$ & $T(1,2,5);41$ & $P_6; 21$\\
\hline
$7$ & $L_{7,6}; 60$ & $T(2,3,2); 49$ & $T(2,2,3); 47$ & $T(1,2,4); 32$ & $P_5; 15$ & no graph exists\\
\hline
$6$ & $K_{1,5}; 37$ & $T(2, 2, 2); 28$ & $T(1,2,3); 24$ & $P_4; 10$ & no graph exists & no graph exists\\
\hline
\end{tabular}
\end{center}
\caption{ The graphs having minimum number of connected subgraphs over $\mathfrak{C}_{n,k}(k)$, $n\leq 12$} \label{Min_n_12}\label{Table}
\end{table}

If $G\in \mathfrak{C}_{n,k}$, $1\leq n\leq 5$, then $k$ is at most $3$, and the graphs with the minimum number of connected subgraphs in $\mathfrak{C}_{n,k}$ are trees. Further, the minimal tree can be characterised using Proposition \ref{min on trees}. For $6\leq n\leq  12$ and various values of $k$, the graphs having the minimum number of connected subgraphs are listed in Table \ref{Table}. The entry ($G ;$ number) in the cells, represent (the graph $G$ attaining the minimum; $F(G)$), respectively, for the corresponding $n$ and $k$. Note that $\mathfrak{C}_{n,k}(k)=\mathfrak{C}_{n,k}$ for $k={1,2,3}$. $\mathfrak{C}_{n,4}$ may contain a graph with girth $3$. So $\mathfrak{C}_{n,4}\neq \mathfrak{C}_{n,4}(4)$. However, by Lemma \ref{triangle free_F}, a graph of girth $3$ can not have the minimum number of connected subgraphs in $\mathfrak{C}_{n,4}$, implying that if a graph has minimum number of connected subgraphs in $C_{n,4}$, then it must belong to $\mathfrak{C}_{n,4}(4)$. Hence Theorem \ref{girth_k}, along with Table \ref{Table}  characterises the minimal graphs in $C_{n,k}$, $1\leq k\leq 4$.

\section{Conclusions}\label{Conclusion}
 We studied the problem of characterising graphs with the minimum number of connected subgraphs within the family of graphs on $n$ vertices and $k$ cut vertices. We solved this problem for $k\leq 4$, showing that for $n\geq 13$, the graph $L_{n,n-k}$ attains the minimum and there is a closed-form formula  (\ref{F(L(n,n-k))}) to compute this minimum. For $n\leq 12$, the minimum value and the graphs attaining the minimum are presented in Table \ref{Table}.  For $k\geq 5$, the minimal graph is characterised in a subfamily of graphs with girth at least $k$.  However, the characterisation of graphs attaining the minimum number of connected subgraphs in the original family $\mathfrak{C}_{n,k}$ for $k\geq 5$ remains open and appears to be a promising direction for future research. In that case, there could be more than one graph attaining the minimum, and so the technique of proof by induction used in this paper may not work. \\

\noindent{\bf Acknowledgements:}  The first author’s research was funded by the Einwechter Centre for Supply Chain Management (a research centre funded by Mr. Dan Einwechter) at Wilfrid Laurier University.

\end{document}